\documentclass[11pt]{amsart}
\usepackage{amsmath,amssymb,txfonts}
\usepackage{amssymb}
\usepackage{amsmath}
\usepackage{mathrsfs}
\usepackage[shortlabels]{enumitem}
\usepackage{varioref}
\usepackage{float}
\usepackage{color}
\usepackage{extarrows}
\def \nin{\notin}
\newtheorem{theorem}{Theorem}[section]
\newtheorem{claim}{Claim}

\theoremstyle{definition}

\newtheorem{problem}[theorem]{Problem}
\theoremstyle{remark}

\numberwithin{equation}{section}

\begin{document}

\title[Sum of Four Cubes]{Problems related to  Waring$-$Goldbach problem involving cubes of primes}

\author[Z. Zhai]{Zhichun Zhai}
\address[Z. Zhai]{Department of Mathematics and Statistics,  MacEwan University,   Edmonton, Alberta T5J2P2 Canada
}
\email{zhaiz2@macewan.ca}

\subjclass[2000]{11P32, 11P55, 11N36} 
\date{}

\dedicatory{}

\keywords{Waring–Goldbach problems, Cubes of primes}

\begin{abstract}
In this note, we try to understand the recent development on the Waring-Goldbach problem involving cubes of primes. Especially, we want to determine whether integers that are either primes,  squares of primes,  cubes of primes,  or a cube of an even number can be written as the sum of four cubes of primes.  Meanwhile, we raise some problems that may deepen our understanding of the problem about the sum of four cubes of primes.   Moreover, some examples suggest that almost all the cubes of integers can be written as the sum of cubes of four integers.   
\end{abstract}

\maketitle

\tableofcontents \pagenumbering{arabic}



\section{Introduction}
Waring's problem, which is closely related to the Fermat Last Theorem \cite{TaylorWiles,Wiles}, is one of the most famous problems is number theory. 
It asks one to determine for each integer $k \geq 2,$ the smallest
positive integer $G(k)$ such that every sufficiently large positive integer is the sum of at most
$G(k)$ many $k-$th powers of positive integers. The best known result is $G(3) \leq 7$ for the cubes,
as established by Linnik \cite{Linnik}. While,  it is  conjectured  that  all  sufficiently  large  integers are the sum of at most four cubes of primes, i.e.
$$N=p_1^3+p_2^3+p_1^3+p_4^3,$$
that  is $G(3)=4.$ 
Such a strong conjecture is out of reach at present.  A theorem  of  Hua  \cite{Hua1, Hua2}   states  that  almost  all  positive  integers satisfying  some necessary  congruence  conditions  are  the  sum  of  five  cubes  of  primes,  
 see e.g. Br\"{u}dern \cite{Brudern}.  Ren \cite{Ren1, Ren2} proved   that the conjecture for sum of four cubes is true for  at least 1.5\% of the positive integers satisfying the necessary conditions. This percentage  was improved by Liu \cite{Liu} to 2.911\%.
Recently, this percentage was proved by Ching and  Tsang \cite{ChineTsang} to be at least 8.25\%. Their result reads as.

\begin{theorem}\cite[Theorem 2]{ChineTsang}\label{Them1}
Let $q_0=3\prod_{p<M} p$ where $M=1000.$ Let $r\geq0$ be an integer such that 
$$2|r,\quad r\not\equiv \pm 1, \pm 3 \ (mod\ 9)\quad \hbox{and}\quad r\not\equiv \pm 1\ (mod\ 7).$$
Then there are at least 8.25\% of the positive integers $N\equiv r\ (mod\ q_0)$ which are sums of  4 cubes of primes.
\end{theorem}

In this note, we would like to understand the type of integers $N$ which can be written as the sum of four cubes of primes.   It is clear that both $q_0$ and $r$ are even. Thus, in  \cite[Theorem 2]{ChineTsang}, integer $N$ are only even.
However, we firstly observe that $N$ can not be some odd numbers such as the cube of primes. Then,  we consider  $N$ that could be square of a prime or $N$ could be a prime. Finally, we will consider the problem about the sum of four cubes of integers.

\section{Sum of Four Cubes of  Primes}
There exist primes  $p,p_1,p_2,p_3$ and $p_4$  satisfy
\begin{equation}\label{eq1}
    p=p_1^3+p_2^3+p_3^3+p_4^3
\end{equation}
or 
\begin{equation}\label{eq2}
    p^2=p_1^3+p_2^3+p_3^3+p_4^3.
\end{equation}

Some primes satisfying (\ref{eq1}) are listed in Table \ref{table1}. 
Some primes satisfying (\ref{eq2}) are listed in Table \ref{table2}. 

\begin{claim}\label{prop1}
There do not exist primes $p,p_1,p_2,p_3,p_4$ satisfying the equation 
\begin{equation}\label{eq3}
    p^3=p_1^3+p_2^3+p_3^3+p_4^3.
\end{equation} 
\end{claim}
Note that  one of $p_1,p_2,p_3,p_4$ must be the smallest prime 2. Without loss of generality, we assume $p_4=2.$ 
So,  if we can conclude that 
$$p^3-2^3=(p-2)(p^2+2p+4)=p_1^3+p_2^3+p_3^3$$ which  holds only when  $p_1=p_2=p_3,$ then $(p-2)(p^2+2p+4)=3p_1^3$ which implies $p=5$ and $p_1^3=39$ which is impossible.  This would help us to prove Claim \ref{prop1}.

\begin{claim}\label{prop2}
  Primes $p,p_1,p_2,p_3,p_4$ satisfy 
\begin{equation}\label{eq4}
    p^5=p_1^3+p_2^3+p_3^3+p_4^3
\end{equation} 
if and only if $p=p_1=p_2=p_3=p_4=2.$
\end{claim}
We may assume that $p_4=2.$
If we can prove that $p^5=p_1^3+p_2^3+p_3^3+2^3$ holds only when $p=p_1=p_2=p_3=2,$ then the Claim \ref{prop2} is true.  

Claims \ref{prop1}-\ref{prop2}  suggest that  if $N$ is sufficiently large and is a cube/fifth power of   a prime number, then $N$ can not be written as the sum of four cubes of primes.   This would provide counterexamples to  the conjecture  $G(3)=4$ which  means that  all  sufficiently large integers  $N$  can be the sum of four cubes of primes.

Claims \ref{prop1}-\ref{prop2}  also suggest that if $N$ is $m-$ power of a prime number when  $m\in\{6,7,8,\cdots\}$ then $N$ can not be written as the sum of four cubes of primes. Thus, we have the following claim. 

\begin{claim} There are no primes $p,p_1,p_2,p_3$ and $p_4$ satisfy 
\begin{equation}\label{eq5}
    p^m=p_1^3+p_2^3+p_3^3+p_4^3.
\end{equation}with integer $m\geq 6.$
\end{claim}

\begin{table}[H]
\begin{tabular}{|c| c| c| c|c|||c| c| c| c| c|||c| c| c| c|c|} 
 \hline
 $p_1$ & $p_2$ & $p_3$ & $p_4$ & $p$ &  $p_1$ & $p_2$ & $p_3$ & $p_4$ & $p$& $p_1$ & $p_2$ & $p_3$ & $p_4$ & $p$ \\ [0.5ex] 
 \hline\hline
2  &3&3&3&89 &2&2&2&5&149&
2  &2&2&7&367\\
  \hline 2&5&5&5&383& 
 2 &3&5&7&503 &2&5&5&7&601\\
  \hline
 2 &3&7&11&1709 &2&2&2&13&2221&
 2 &3&5&13&2357 \\
  \hline
  2&11&11&11&4001&
 2 &2&2&17&4937 &2&5&5&17&5171\\
  \hline
  2&13&13&13&6599 &2&2&2&19&6883&
 2 &3&5&19&7019 \\
  \hline2&3&7&19&7237&
 2 &5&13&17&7243 &2&11&11&17&7583\\
  \hline
 2 &3&13&19&9091 &2&7&17&17&10177&
  2&13&13&19&11261 \\
  \hline2&3&17&19&11807&
  2&17&17&17&14747 &2&13&19&19&16693\\
  \hline
 2 &7&17&23&17431 &2&2&2&29&24413&
  2&3&7&29&24767 \\
  \hline
  2&11&23&23&25673&
 2 &11&11&29&27059 &2&3&7&31&30169
 \\
  \hline
 2 &11&19&29&32587 &2&3&17&31&34739&
 2 &19&19&31&43517 
 \\
  \hline
  2&23&23&29&48731&
 2 &3&7&37&51031 &2&3&7&37&51347
 \\\hline
 2&7&13&37& 53201&2&11&11&37&53323&
2 &17&29&29&53699
\\
  \hline
  2&23&23&31&54133&
 2&3&31&31&59617&2&23&29&29&60953\\
  \hline
 2&19&31&31&66449&2&17&23&37&67741
&  2&3&13&41&71153
\\
  \hline
  2&5&19&41&76039
& 2&29&29&31&78577&2&2&2&43&79531
  \\
  \hline
 2&17&19&41&80701&2&7&23&41&81439
 &2&7&23&41&81439\\
  \hline2&31&31&31&89381
  & 2&5&23&43&91807&2&23&23&41&93263
  \\
  \hline
  2&11&29&41&94649&2&29&29&37&99439
  & 2&3&37&37&101341
  \\
  \hline2&13&37&37&103511
  & 2&5&17&47&113657
  &2&17&17&47&108869
  \\
  \hline
 2&17&17&47&115603&2&7&23&47&116341
 &2 &11&23&47&117329\\
  \hline
\hline
\end{tabular}
\caption{\label{table1}
Examples of prime solutions of (\ref{eq1}) with $b\leq 120000$}
\end{table}

\begin{table}[H]
\begin{tabular}{|c| c| c| c| c|||c| c| c| c| c|||c| c| c| c| c|}
 \hline
 $p_1$ & $p_2$ & $p_3$ & $p_4$ & $p$& $p_1$ & $p_2$ & $p_3$ & $p_4$ & $p$&$p_1$ & $p_2$ & $p_3$ & $p_4$ & $p$\\ [0.5ex] 
 \hline\hline
2&29  &47  &109    &1193 &
2&13  &107  &137    &1949 &
2&79  &107  &311    &5639 
\\
\hline
2&11  &313  &317    &7907 &
 2&173&233&379&85012& 2& 79  &197  &461    &10301\\
\hline
2&317&547&593&20101 &2&107&239&751&20939
&
2&439&617&809&29137
\\
\hline
2&23&107&967&30091
&2&107&499&1019&34403&2&239&977&1489&65173
\\
\hline
2&317&743&1621&68567&2&373&1319&1373&70249
&2&317&701&1663&70537\\
\hline
2&919&1187&1511&76801
&2&491&1481&1549&84163&2&317&1213&2111&105943
\\
\hline
2&227&1481&2053&109147&2&547&743&2351&116483
&2&1481&2081&2239&153247
\\
\hline
2&499&569&2957&161753&2&83&1787&2851&169943&2&599&977&3049&171733
\\
\hline
2&1123&1493&2969&180563&2& 71&1871&3331&208589
&2&1657&2003&3221&214483
\\
\hline2&1499&2579&3067&222197
&2&1367&2621&3109&224969&2&449&631&3881&242483
\\
\hline
2&2137&2699&3389&261427&2&2203&3203&3803&313933
&2&1123&1493&2969&175829
\\
\hline
2 &1471&1481&2969&180563&
2&71&1871&3331&208589&2&1657&2003&3221&214483\\
\hline
2&1499&2579&3067&222197&2&1367&2621&3109&224969
&2&449&631&3881&242483
\\
\hline
2&2137&2699&3389&261427&
2&509&1877&4159&280507&2&2203&3203&3803&313933
\\
\hline
2&1627&2153&4421&317327&2&947&3719&4231&357809
&2&317&653&5167&371831
\\
\hline2&2927&4283&4327&429719
&2&653&4127&5101&450887&2&2903&3851&5527&500413
\\
\hline
2&3167&4337&5471&526459&2&281&4463&5857&538367
&2&2633&5167&5417&561389
\\
\hline2&2087&5171&5623&570217 &&&&&&&&&&
\\
\hline
\end{tabular}
\caption{\label{table2}Examples of prime solutions of (\ref{eq2}) with $p_i\leq  6000, i=1,2,3,4.$}

\end{table}

Our next question is stated as follows. 

\begin{problem} Among all $N$ that can be written as the sum of four cubes of primes in Theorem \ref{Them1},  how many of them are cubes of an integer? 
\end{problem}
This problem asks for the number of  integers $b$ such that
\begin{equation}\label{eq6}
b^3=p_1^3+p_2^3+p_3^3+p_4^3
\end{equation}
holds for primes $p_1,p_2,p_3,p_4.$ 
Some solutions of equation (\ref{eq6})  with positive  integer $b$ and primes $p_1,p_2,p_3,p_4\leq 3000$ are listed in   Table (\ref{tab3}). Table (\ref{tab3}) suggests that among all $N$ that can be written as the sum of four cubes of primes, only a very small portion of them are the cube of primes.

In \cite{Bruedern},  the set $\mathcal{N}$ was defined as the set of all natural numbers n satisfying
$$2|n,\quad n\not\equiv \pm 1, \pm 3 \ (mod\ 9)\quad \hbox{and}\quad r\not\equiv \pm 1\ (mod\ 7).$$
Note that if a number $n\nin \mathcal{N} $ can be written as a sum of four cubes of primes, then one
of the cubes in the representation must be that of 2, 3 or 7.

\begin{claim}\label{claim2}
If $n$ is can 
be written as a sum of four cubes of primes with one prime $2,3,$ or $7,$ then  $n\nin \mathcal{N} .$
\end{claim}

The examples in Table \ref{tab3} suggest that 
Claim  \ref{claim2} is true. But, we do not have a rigorous proof.

\begin{table}[H]
\begin{tabular}{|c| c| c| c| c|||c| c| c| c| c|} 
 \hline
 $p_1$ & $p_2$ & $p_3$ & $p_4$ & $b$& $p_1$ & $p_2$ & $p_3$ & $p_4$ & $b$\\ [0.5ex] 
 \hline\hline
 3 &3  &7  &11  &12&
 3  &317   & 379  &509  &602\\ 
 \hline
5  &67  &191   &313 &336&
7  &29  &967   &1013 &1248 \\\hline
 7&43 & 47   &47 &66 &
47  &281  &283  &541 &588 \\  \hline
53  &307  &389   &1627 &1638 &
59  &179  &421   &1069 &1092 \\\hline
 107  &347  &349  &349 &504 &
181  &463  & 491  &593 &756 \\
\hline
337  &1061  &1063   &1571 &1848 &
53  &307  &389   &1627 &1638 \\
\hline
353  &617  &1123   &2011 &2142 &
661  &887  &1033   &2027 & 2184\\
  \hline
 2 &433  &797   &2251 &2289 &
101  &281  &397   &2389 &2394 \\
  \hline
43  &857  &1861   & 2423&2772 & &&&&\\
\hline
\end{tabular}
\caption{\label{tab3}
Solutions to $b^3=p_1^3+p_2^3+p_3^3+p_4^3$
with positive integers $b$ and  primes less than 3000.}
\end{table}

\section{Sum of  Four Cubes of Integers}
 Hardy and Littlewood in \cite{Hardy}  proved that almost  positive integers
are representable as the sum of five positive integral cubes. 
  Davenport in \cite{Davenport} 
proved that not more than $0.134\cdots$ of N positive integers less than N are
representable as the sum of three positive integral cubes.  Instead,   we consider how many cubes can be written as the sum of four cubes of integers. 

\begin{problem} Among all integers that can be written as the sum of four cubes of positive integers,  how many of them are the cube of an integer? 
\end{problem}
This problem asks for the number of  integers $b$ such that
\begin{equation}\label{eq7}
b^3=x_1^3+x_2^3+x_3^3+x_4^3
\end{equation}
holds for integers $b, x_1,x_2,x_3,x_4.$ 
Some solutions of equation (\ref{eq7})  with positive  integers $b,$  $x_1,x_2,x_3,x_4\leq 1000$ are listed in   the following table.

\begin{table}[H]
\begin{tabular}{|c| c| c| c| c|||c| c| c| c| c|||c| c| c| c| c|} 
 \hline
 $x_1$ & $x_2$ & $x_3$ & $x_4$ & $b$ & $x_1$ & $x_2$ & $x_3$ & $x_4$ & $b$& $x_1$ & $x_2$ & $x_3$ & $x_4$ & $b$\\ [0.5ex] 
 \hline\hline
 1 & 1 & 5 &6  &7& 3&3 &7 &11 &12 & 
1 &5 &7 &12 &13\\  \hline
5 &7 &9 &10 &13&
 2& 2&10 &12 &14 & 
1 &7 &614 &14 &18\\
 \hline
 4& 7&8 &17 &18& 
 11& 12& 13&14 &20&
3 &3 &15 &18 &21\\
 \hline
 6& 14&16 &19 &24&
7 &9 &9 &24 &26&
2 &10 &14 & 22&26\\
 \hline
3 &15 & 17&21 &26&
10 &14 & 18&20 &26&
1&10&15  &26 &28\\
 \hline
2  &4 &13 &27 &28&
4 &4 &20 &24 &28&
4 &6 &16 &26 &28\\
 \hline
9 &13 &19 &23 &28&
13& 15&21 &23 &30 &
9& 13&17 &28 &31\\
 \hline
7&8 &17 &30 &32&
3&14 &15 &31 &33&
11& 13&24 &28 &34\\
 \hline
5&5 &25 &30 &35&
8&14 &16 &34 &36&
9&9&21  &33 &36\\
 \hline
6& 16& 28& 29&37&
3& 5&28 &32 &38&
2& 7& 16&38 &39\\
 \hline
3&15 &21 &36 &39&
7
&9 &15 &38 &39&
15&21 &27 &30 &39\\
 \hline
22&24 &26 &28 &40&
4&12 &15 &41 &42&
6&6 &30 &36 &42 \\
 \hline
 6&9 &24 &39 &42&
6&15& 25 &38 &42 & 17&21 &21 &37 &942\\
 \hline
15&28 &29 &31 &43&
16&24&29  &35 &44
&
10& 13&14 &44 &45\\
 \hline
12 &28&30&36&46&
7 &17&20&46&48&
10 &28&32&38&48
 \\
 \hline
13 &22&35&38&48&
20 &21&24&47&49&
1 &1&24&47&49
 \\
 \hline
 2&13&33&43&49&
7& 7&25&42&49&
19& 23&34&39&49
 \\
 \hline
 14&18&18&48&50&
14 &30&33&39&50&
7 &16&29&47 &51
 \\
 \hline
 18&27&28&44&51&
4 &20&28&48&52&
6 &30&34&42&52
 \\
 \hline
 10&22&38&42&52&
19&21&38&42&52&
20 &28&36&40&52
 \\
 \hline
 3&35&36&39&53&
27&30&35&39&53&
3 &21&42&42&54
 \\
 \hline
 7&17&38&46&54&
 12&21&24&51&54&
21& 24&38&43&54
 \\
 \hline
23 &24&24&49&54&
27&29&34&42&54&
7 &12&34&50&55
 \\
 \hline
7&24&38&46&55&
7 &19&28&51&55&
2 &20&30&52&56
 \\
 \hline
4 &8&26&54&56&
8&8&40&48&56&
8 &9&20&55&56
 \\
 \hline
8 &12&32&52&56&
10 &26&29&51&56&
18 &26&38&46&56
 \\
 \hline
 20&35&37&42&56&
1 &8&42&48&57&
13 &23&28&53&57
 \\
 \hline
13 &21&41&50&59&
14 &19&44&48&59&
21& 23&26&55&59
 \\
 \hline

 10&23&44&49&60&
15 &15&35&55&60&
26 &30&42&46&60
 \\
 \hline
 33&36&39&42&60&
5 &6&35&58&62&
18 &26&34&56&62
 \\
 \hline
 2&15&27&61&63&
3 &28&35&57&63&
5 &9&45&54&63
 \\
 \hline
 9&9&45&54&63&
12 &39&40&50&63&
13 &23&27&60&63
 \\
 \hline
31 &33&39&50&63&
1 &16&20&63&64&
8 &15&49&52&64
 \\
 \hline
 14&16&34&60&64&
 21&32&43&56&65&
 1&24&44&56&65
  \\
 \hline
 4&5&29&63&65&
5 &25&35&60&65&
8 &17&42&58&65
  \\
 \hline
 11&21&41&58&65&
25 &35&45&50&65&
 4&26&35&61&66
  \\
 \hline
 6&23&49&54&66&
 6&28&30&62&66&
 7&43&47&47&66
  \\
 \hline
 8&10&33&63&66&
 13&16&35&62&66&
 13&18&42&59&66
  \\
 \hline
 20&24&25&63&66&
 31&38&44&49&66&
3 &22&42&60&67
  \\
 \hline
 7&17&43&60&67&
22 &22&42&59&68&
13 &32&42&59&68
  \\
 \hline
 19&35&39&59&68&
22 &26&48&56&68&
 26&30&35&61&68
  \\
 \hline
 6&28&53&54&69&
11 &41&49&52&69&
18 &42&45&54&69
  \\
 \hline
 20&35&37&62&69&
26 &28&37&62&69&
10 &10&50&60&70
  \\
 \hline
10 &15&50&65&70&
 11&16&48&65&70&
13 &19&39&65& 70
  \\
 \hline
 15&24&52&57&70&
16&22&24&68&70&
19 &36&52&53&70
  \\
 \hline
 25&39&48&54&70&
 3&18&56&56&72&
4 &28&56&56&72
  \\
 \hline

\end{tabular}
\caption{\label{table4}
Positive integers solutions of equation (\ref{eq7}) with $b\leq 72$}
\end{table}

We checked all  $b\leq 1000,$  we observe that when  $b$ takes the following values $$1,2,3,4,5,6,8,9,10,11,15,16,17,19,22, 27,29, 47,58,  61,71,$$ then $b^3$ can not be represented as the sum of four cubes of positive integers.   Table \ref{table4} lists solutions to (\ref{eq7}) with $b\leq 72.$ 
We propose  the following claim.

\begin{claim}
For almost all sufficiently large number $b,$ $b^3$ can be represented as the sum of four cubes of positive integers. The representation is not unique. 
\end{claim}


\begin{thebibliography}{30}

\bibitem{Brudern}
J. Br \"{u}dern, A sieve approach to the Waring–Goldbach problem (I): sums of four cubes, Ann. Sci.E \'{c}ole Norm. Sup. 28 (1995) 461–476.

\bibitem{Bruedern}
J. Br \"{u}dern,  and K. Kawada. On the Waring–Goldbach problem for cubes. Glasgow Mathematical Journal 51.3 (2009): 703-712.


\bibitem{ChineTsang}
T. W. Ching,  and K. M. Tsang, Waring–Goldbach problem involving cubes of primes. Mathematische Zeitschrift 297.3 (2021): 1105-1117.

\bibitem{Davenport}
  H. Davenport, On Waring’s problem for cubes, Acta Math. 71 (1939) 123–143.
  
  \bibitem{Hardy}
  G H Hardy and J E Littlewood, 
  Some problems of 'Partitio numerorum' (VI): Further researches in Waring's Problem, Mathematische Zeitschrift 23 (1925), 1-37.


  \bibitem{Hua1} 
  L.K.  Hua,  Some  results  in  additive  prime  number  theory,  Quart.  J.  Math.  (Oxford)  9  (1938)68–80.
  
  \bibitem{Hua2}
  L.K.  Hua, Additive  Theory of  Prime  Numbers, Science  Press,  Beijing,  1957 (in Chinese);  Englishversion, Amer. Math. Soc., Providence, RI, 1965

\bibitem{Linnik}
Y. Linnik,  On the representation of large numbers as sums of seven cubes. Mat. Sb. 12, 218–224 (1943)


\bibitem{Liu}
Z. Liu, Density of the sums of four cubes of primes. Journal of Number Theory 132.4 (2012): 735-747.

\bibitem{Ren1}
 X.M. Ren, Density of integers that are the sum of four cubes of primes, Chin. Ann. Math. Ser. B 22 (2001) 233–242.

\bibitem{Ren2} 
X.M. Ren, Sums of four cubes of primes, J. Number Theory 98 (2003) 156–171.

\bibitem{TaylorWiles}

R. Taylor, and A.  Wiles,   Ring-theoretic properties of certain Hecke algebras. Annals of Mathematics, (1995). 553-572.


\bibitem{Wiles}
A. Wiles,  Modular elliptic curves and Fermat's last theorem. Annals of mathematics 141.3 (1995): 443-551.


\end{thebibliography}
\end{document}